\documentclass[a4paper,12pt]{article}
\usepackage[all]{xy}
\usepackage[T1]{fontenc}
\usepackage[dvips]{graphicx}
\usepackage{amsmath,amssymb,amsfonts,amsthm,latexsym,mathrsfs,textcomp,verbatim}
\xyoption{web}

\title{Yoneda representations of flat functors and classifying toposes}
\author{{Olivia Caramello} \vspace{3 mm}\\ {\small DPMMS, University of Cambridge,}\\{\small Wilberforce Road, Cambridge CB3 0WB, UK}\\{\small O.Caramello@dpmms.cam.ac.uk}}
\date{\today}
\begin{document}
\bgroup           
\let\footnoterule\relax  
\maketitle
\flushleft  
\begin{abstract}
In this paper, we first introduce a technique that we call ``Yoneda representation of flat functors'', based on ideas from indexed category theory; then we provide applications of this technique to the theory of classifying toposes. Specifically, we obtain results characterizing the models of a theory classified by a topos of the form ${\bf Sh}({\cal C},J)$ in terms of the models of a theory classified by the topos $[{\cal C}^{\textrm{op}},{\bf Set}]$.
\end{abstract} 
\egroup 
\flushleft
\vspace{5 mm}


\def\Monthnameof#1{\ifcase#1\or
   January\or February\or March\or April\or May\or June\or
   July\or August\or September\or October\or November\or December\fi}
\def\today{\number\day~\Monthnameof\month~\number\year}

%
%
%
\def\pushright#1{{
   \parfillskip=0pt            
   \widowpenalty=10000         
   \displaywidowpenalty=10000  
   \finalhyphendemerits=0      
  %
   \leavevmode                 
   \unskip                     
   \nobreak                    
   \hfil                       
   \penalty50                  
   \hskip.2em                  
   \null                       
   \hfill                      
   {#1}                        
  %
   \par}}                      

\def\qed{\pushright{$\square$}\penalty-700 \smallskip}

\newtheorem{theorem}{Theorem}[section]

\newtheorem{proposition}[theorem]{Proposition}

\newtheorem{scholium}[theorem]{Scholium}

\newtheorem{lemma}[theorem]{Lemma}

\newtheorem{corollary}[theorem]{Corollary}

\newtheorem{conjecture}[theorem]{Conjecture}

\newenvironment{proofs}%
 {\begin{trivlist}\item[]{\bf Proof }}%
 {\qed\end{trivlist}}

  \newtheorem{rmk}[theorem]{Remark}
\newenvironment{remark}{\begin{rmk}\em}{\end{rmk}}

  \newtheorem{rmks}[theorem]{Remarks}
\newenvironment{remarks}{\begin{rmks}\em}{\end{rmks}}

  \newtheorem{defn}[theorem]{Definition}
\newenvironment{definition}{\begin{defn}\em}{\end{defn}}

  \newtheorem{eg}[theorem]{Example}
\newenvironment{example}{\begin{eg}\em}{\end{eg}}

  \newtheorem{egs}[theorem]{Examples}
\newenvironment{examples}{\begin{egs}\em}{\end{egs}}


\mathcode`\<="4268  
\mathcode`\>="5269  
\mathcode`\.="313A  
\mathchardef\semicolon="603B 
\mathchardef\gt="313E
\mathchardef\lt="313C

\newcommand{\app}
 {{\sf app}}

\newcommand{\Ass}
 {{\bf Ass}}

\newcommand{\ASS}
 {{\mathbb A}{\sf ss}}

\newcommand{\Bb}
{\mathbb}

\newcommand{\biimp}
 {\!\Leftrightarrow\!}

\newcommand{\bim}
 {\rightarrowtail\kern-1em\twoheadrightarrow}

\newcommand{\bjg}
 {\mathrel{{\dashv}\,{\vdash}}}

\newcommand{\bstp}[3]
 {\mbox{$#1\! : #2 \bim #3$}}

\newcommand{\cat}
 {\!\mbox{\t{\ }}}

\newcommand{\cinf}
 {C^{\infty}}

\newcommand{\cinfrg}
 {\cinf\hy{\bf Rng}}

\newcommand{\cocomma}[2]
 {\mbox{$(#1\!\uparrow\!#2)$}}

\newcommand{\cod}
 {{\rm cod}}

\newcommand{\comma}[2]
 {\mbox{$(#1\!\downarrow\!#2)$}}

\newcommand{\comp}
 {\circ}

\newcommand{\cons}
 {{\sf cons}}

\newcommand{\Cont}
 {{\bf Cont}}

\newcommand{\ContE}
 {{\bf Cont}_{\cal E}}

\newcommand{\ContS}
 {{\bf Cont}_{\cal S}}

\newcommand{\cover}
 {-\!\!\triangleright\,}

\newcommand{\cstp}[3]
 {\mbox{$#1\! : #2 \cover #3$}}

\newcommand{\Dec}
 {{\rm Dec}}

\newcommand{\DEC}
 {{\mathbb D}{\sf ec}}

\newcommand{\den}[1]
 {[\![#1]\!]}

\newcommand{\Desc}
 {{\bf Desc}}

\newcommand{\dom}
 {{\rm dom}}

\newcommand{\Eff}
 {{\bf Eff}}

\newcommand{\EFF}
 {{\mathbb E}{\sf ff}}

\newcommand{\empstg}
 {[\,]}

\newcommand{\epi}
 {\twoheadrightarrow}

\newcommand{\estp}[3]
 {\mbox{$#1 \! : #2 \epi #3$}}

\newcommand{\ev}
 {{\rm ev}}

\newcommand{\Ext}
 {{\rm Ext}}

\newcommand{\fr}
 {\sf}

\newcommand{\fst}
 {{\sf fst}}

\newcommand{\fun}[2]
 {\mbox{$[#1\!\to\!#2]$}}

\newcommand{\funs}[2]
 {[#1\!\to\!#2]}

\newcommand{\Gl}
 {{\bf Gl}}

\newcommand{\hash}
 {\,\#\,}

\newcommand{\hy}
 {\mbox{-}}

\newcommand{\im}
 {{\rm im}}

\newcommand{\imp}
 {\!\Rightarrow\!}

\newcommand{\Ind}[1]
 {{\rm Ind}\hy #1}

\newcommand{\iten}[1]
{\item[{\rm (#1)}]}

\newcommand{\iter}
 {{\sf iter}}

\newcommand{\Kalg}
 {K\hy{\bf Alg}}

\newcommand{\llim}
 {{\mbox{$\lower.95ex\hbox{{\rm lim}}$}\atop{\scriptstyle
{\leftarrow}}}{}}

\newcommand{\llimd}
 {\lower0.37ex\hbox{$\pile{\lim \\ {\scriptstyle
\leftarrow}}$}{}}

\newcommand{\Mf}
 {{\bf Mf}}

\newcommand{\Mod}
 {{\bf Mod}}

\newcommand{\MOD}
{{\mathbb M}{\sf od}}

\newcommand{\mono}
 {\rightarrowtail}

\newcommand{\mor}
 {{\rm mor}}

\newcommand{\mstp}[3]
 {\mbox{$#1\! : #2 \mono #3$}}

\newcommand{\Mu}
 {{\rm M}}

\newcommand{\name}[1]
 {\mbox{$\ulcorner #1 \urcorner$}}

\newcommand{\names}[1]
 {\mbox{$\ulcorner$} #1 \mbox{$\urcorner$}}

\newcommand{\nml}
 {\triangleleft}

\newcommand{\ob}
 {{\rm ob}}

\newcommand{\op}
 {^{\rm op}}

\newcommand{\pepi}
 {\rightharpoondown\kern-0.9em\rightharpoondown}

\newcommand{\pmap}
 {\rightharpoondown}

\newcommand{\Pos}
 {{\bf Pos}}

\newcommand{\prarr}
 {\rightrightarrows}

\newcommand{\princfil}[1]
 {\mbox{$\uparrow\!(#1)$}}

\newcommand{\princid}[1]
 {\mbox{$\downarrow\!(#1)$}}

\newcommand{\prstp}[3]
 {\mbox{$#1\! : #2 \prarr #3$}}

\newcommand{\pstp}[3]
 {\mbox{$#1\! : #2 \pmap #3$}}

\newcommand{\relarr}
 {\looparrowright}

\newcommand{\rlim}
 {{\mbox{$\lower.95ex\hbox{{\rm lim}}$}\atop{\scriptstyle
{\rightarrow}}}{}}

\newcommand{\rlimd}
 {\lower0.37ex\hbox{$\pile{\lim \\ {\scriptstyle
\rightarrow}}$}{}}

\newcommand{\rstp}[3]
 {\mbox{$#1\! : #2 \relarr #3$}}

\newcommand{\scn}
 {{\bf scn}}

\newcommand{\scnS}
 {{\bf scn}_{\cal S}}

\newcommand{\semid}
 {\rtimes}

\newcommand{\Sep}
 {{\bf Sep}}

\newcommand{\sep}
 {{\bf sep}}

\newcommand{\Set}
 {{\bf Set }}

\newcommand{\Sh}
 {{\bf Sh}}

\newcommand{\ShE}
 {{\bf Sh}_{\cal E}}

\newcommand{\ShS}
 {{\bf Sh}_{\cal S}}

\newcommand{\sh}
 {{\bf sh}}

\newcommand{\Simp}
 {{\bf \Delta}}

\newcommand{\snd}
 {{\sf snd}}

\newcommand{\stg}[1]
 {\vec{#1}}

\newcommand{\stp}[3]
 {\mbox{$#1\! : #2 \to #3$}}

\newcommand{\Sub}
 {{\rm Sub}}

\newcommand{\SUB}
 {{\mathbb S}{\sf ub}}

\newcommand{\tbel}
 {\prec\!\prec}

\newcommand{\tic}[2]
 {\mbox{$#1\!.\!#2$}}

\newcommand{\tp}
 {\!:}

\newcommand{\tps}
 {:}

\newcommand{\tsub}
 {\pile{\lower0.5ex\hbox{.} \\ -}}

\newcommand{\wavy}
 {\leadsto}

\newcommand{\wavydown}
 {\,{\mbox{\raise.2ex\hbox{\hbox{$\wr$}
\kern-.73em{\lower.5ex\hbox{$\scriptstyle{\vee}$}}}}}\,}

\newcommand{\wbel}
 {\lt\!\lt}

\newcommand{\wstp}[3]
 {\mbox{$#1\!: #2 \wavy #3$}}

\newcommand{\fu}[2]
{[#1,#2]}

\newcommand{\st}[2]
 {\mbox{$#1 \to #2$}}

\section{Preliminary facts}
In this section we introduce the terminology and recall the facts from the theory of indexed categories that will be useful for our analysis. 
We refer the reader to \cite{El} (especially sections B1.2, B2.3 and B3.1) and to \cite{PS} for the background.\\

By a topos (defined) over \Set we mean an elementary topos $\cal{E}$ such that there exists a (necessarily unique up to isomorphism) geometric morphism $\stp{\gamma_{\cal{E}}}{\cal{E}}{\Set}$; we denote by $\gamma_{\cal{E}}^{\ast}$ the inverse image functor and by $\Gamma_{\cal{E}}$ its right adjoint, that is the ``global sections'' functor. A topos is defined over \Set if and only if it is locally small and has arbitrary set-indexed copowers of $1$; in particular every locally small cocomplete topos (and hence every Grothendieck topos) is defined over $\Set$.\\  

Given a category $\cal{C}$ and a topos $\cal{E}$ defined over $\Set$, we can always internalize $\cal{C}$ into $\cal{E}$ by means of $\gamma_{\cal{E}}^{\ast}$; the resulting internal category in $\cal{E}$ will be denoted by $\mathbb{C}$.\\
Every topos $\cal{E}$ (over $\Set$) gives rise to a $\cal{E}$-indexed category $\mathbb{E}$ obtained by indexing $\cal{E}$ over itself; the inverse image functor $\gamma_{\cal{E}}^{\ast}$ then induces an indexing of $\cal{E}$ over $\Set$, which coincides with the canonical indexing of $\cal{E}$ provided that $\cal{E}$ is cocomplete and locally small.\\

We will generally denote indexed categories by underlined letters, to distinguish them from their underlying categories which will be denoted by the corresponding simple letters; so for example the underlying category of an indexed category $\mathbb D$ will be denoted by $\cal D$; an exception to this rule will be the notation for indexed categories arising as the indexing of a cartesian category over itself: in this case the indexed category corresponding to a cartesian category $\cal S$ will be simply denoted by $\mathbb S$. Also, internal categories will be denoted by letters $\mathbb C,\mathbb D$, etc. and we will not modify their notation when they are considered as indexed categories.\\

For a topos $\cal E$ and an internal category $\mathbb C$ in $\cal E$, we have a $\cal{E}$-indexed category $\underline{\fu{\mathbb{C}}{\cal{E}}}$, whose underlying category is the category $\fu{\mathbb{C}}{\cal{E}}$ of diagrams of shape $\mathbb{C}$ in $\mathbb{E}$ and morphisms between them. $\fu{\mathbb{C}}{\cal{E}}$ is equivalent (naturally in $\cal{E}$) to the category of $\cal{E}$-indexed functors $\st{\mathbb{C}}{\mathbb{E}}$ and indexed natural transformations between them (by Lemma B2.3.13 in \cite{El}) and also, if $\cal{E}$ is cocomplete and locally small, to the category $\fu{\cal{C}}{\cal{E}}$ (by Corollary B2.3.14 in \cite{El}). For this reason, we will restrict our attention to locally small cocomplete toposes; we will occasionally loosely refer to them simply as toposes.\\ 

The equivalence between $\fu{\mathbb{C}}{\cal{E}}$ and $\fu{\cal{C}}{\cal{E}}$ restricts to an equivalence between the full subcategories $\bf{Tors}(\mathbb{C},\cal{E})$ of $\mathbb{C}$-torsors in $\cal{E}$ (as in section B3.2 of \cite{El}) and $\bf Flat(\cal{C}, \cal{E})$ of flat functors $\st{\cal{C}}{\cal{E}}$ (as in chapter VII of \cite{MM}).\\ Given a functor $F\in \fu{\cal{C}}{\cal{E}}$, the internal diagram that corresponds to it via the natural equivalence $\fu{\cal{C}}{\cal{E}}\simeq \fu{\mathbb{C}}{\cal{E}}$ will be called the internalization of $F$ and denoted by $F^{i}$; of course, this is defined only up to isomorphism.\\

The $\cal{E}$-indexed category $\underline{\fu{\mathbb{C}}{\cal{E}}}$ is locally small (by Lemma B2.3.15 in \cite{El}); from this it follows that there exists a $\cal{E}$-indexed hom functor\\ $\stp{ Hom_{\underline{\fu{\mathbb{C}}{\cal{E}}}}^{\cal{E}}}{\underline{\fu{\mathbb{C}}{\cal{E}}}^{\textrm{op}}\times\underline{[\mathbb{C},\cal{E}]}}{\mathbb{E}}$ whose underlying functor\\
$\stp{Hom_{\fu{\mathbb{C}}{\cal{E}}}^{\cal{E}}}{\fu{\mathbb{C}}{\cal{E}}^{\textrm{op}}\times [\mathbb{C},\cal{E}]}{\cal{E}}$
assigns to a pair of diagrams $F$ and $G$ in $\fu{\mathbb{C}}{\cal{E}}$ an object $Hom_{\fu{\mathbb{C}}{\cal{E}}}^{\cal{E}}(F, G)$ of $\cal{E}$, which we call the object of morphisms from $F$ to $G$ in $\underline{\fu{\mathbb{C}}{\cal{E}}}$.
Also, there is a Yoneda $\cal E$-indexed functor $\stp{\mathbb Y}{\mathbb{C}}{\underline{\fu{{\mathbb{C}}^{\textrm{op}}}{\cal{E}}}}$, which plays in this context the same role as that of a Yoneda functor in ordinary category theory.\\ 

We denote by $\stp{E^{\ast}}{\cal{E}}{{\cal E}}\slash E$ the pullback functor along the unique arrow $\st{E}{1}$, that is the (logical) inverse image functor of the local homeomorphism $\st{{\cal{E}}\slash{E}}{\cal{E}}$. Then, from the equivalence $\fu{\cal{C}}{\cal{E}}\simeq \fu{\mathbb{C}}{\cal{E}}$ we deduce the existence of a hom functor $\stp{Hom_{\fu{\cal{C}}{\cal{E}}}^{\cal{E}}}{[\cal{C},\cal{E}]^{\textrm{op}} \times \fu{\cal{C}}{\cal{E}} }{\cal{E}}$, which assigns to each pair of functors $F$ and $G$ in $\fu{\cal{C}}{\cal{E}}$ an object $Hom_{\underline{\fu{\mathbb{C}}{\cal{E}}}}^{\cal{E}}(F,G)$ of $\cal{E}$ such that for each $E\in \cal{E}$ the morphisms $\st{E}{Hom_{\underline{\fu{\mathbb{C}}{\cal{E}}}}^{\cal{E}}(F,G)}$ in $\cal{E}$ are in natural bijection with the morphisms in $[\cal{C},\cal{E}$$ \slash E]$ from $E^{\ast}\circ F$ to $E^{\ast}\circ G$, that is with the natural transformations $E^{\ast}\circ F \imp E^{\ast}\circ G$.\\
We remark that, since $\bf Flat(\cal{C}, \cal{E})$ is a full subcategory of $\fu{\cal{C}}{\cal{E}}$, we may use the objects $Hom_{\underline{\fu{\mathbb{C}}{\cal{E}}}}^{\cal{E}}(F,G)$ for $F, G$ $\in \bf Flat(\cal{C}, \cal{E})$ as the objects of morphisms from $F$ to $G$ in $\underline{\bf Flat(\cal{C}, \cal{E})}$.\\ 

Given a $\cal S$-indexed category $\underline{\mathbb D}$ and an object $I\in S$, we have a ${\cal S}\slash I$-indexed category ${\underline{\mathbb D}}\slash I$ (defined in the obvious way), which is called the localization of $\underline{\mathbb D}$ at $I$. If $\underline{\mathbb D}$ and $\underline{\mathbb E}$ are two $\cal S$-indexed categories, we denote by $[\underline{\mathbb D},\underline{\mathbb E}]$ the category of $\cal S$-indexed functors from $\underline{\mathbb D}$ to $\underline{\mathbb E}$ and indexed natural transformations between them. The assignment $I\rightarrow [{\underline{\mathbb D}}\slash I,{\underline{\mathbb E}}\slash I]$ is pseudofunctorial in $I\in \cal S$ and makes $[\underline{\mathbb D},\underline{\mathbb E}]$ into a $\cal S$-indexed category.

\section{Yoneda representations}
It is well known that, by Yoneda, for each $F \in [{\cal{C}}^{\textrm{op}}, \Set]$ there is a natural isomorphisms of functors\\
\[
F \cong Hom_{\fu{{\cal C}^{\textrm{op}}}{\Set}}^{\Set}(Y(-),F),
\]
where $Hom_{\fu{{\cal{C}}^{\textrm{op}}}{\Set}}^{\Set}(Y(-),F)$ is the functor given by the composite
\[
{\cal{C}}^{\textrm{op}} \stackrel{Y^{\textrm{op}} \times \Delta F}{\longrightarrow} [{\cal C}^{\textrm{op}}, \Set]^{\textrm{op}}\times [{\cal C}^{\textrm{op}}, \Set]
\stackrel{Hom_{\fu{\cal{C}^{\textrm{op}}}{\Set}}^{\Set}}{\longrightarrow} \Set.  
\]
Thanks to the remarks in the last section we are able to generalize this result to the case of functors with values in an arbitary topos. In fact, the following result holds.

\begin{theorem}
Let $\cal{C}$ a small category and $\cal{E}$ be a locally small cocomplete topos. Then for every functor $\stp{F}{{\cal C}^{\textrm{op}}}{\cal E}$, there is a natural isomorphism of functors
\[
F \cong Hom_{\fu{{\cal C}^{\textrm{op}}}{\cal E}}^{\cal E}(\overline{Y}(-),F),
\]
where $\stp{\overline{Y}}{\cal C}{[{\cal C}^{\textrm{op}},{\cal E}]}$ is the functor given by the composite
\[
{\cal{C}} \stackrel{Y}{\longrightarrow} [{\cal C}^{\textrm{op}}, \Set]
\stackrel{\gamma_{\cal{E}}^{\ast} \circ -}{\longrightarrow} [{\cal C}^{\textrm{op}}, {\cal E}] 
\]
and $Hom_{\fu{{\cal C}^{\textrm{op}}}{\cal E}}^{\cal E}(\overline{Y}(-),F)$ is the functor given by the composite
\[
{\cal{C}}^{\textrm{op}} \stackrel{\overline{Y}^{\textrm{op}} \times \Delta F  }{\longrightarrow} [{\cal C}^{\textrm{op}}, {\cal E}]^{\textrm{op}}\times [{\cal C}^{\textrm{op}}, {\cal E}]
\stackrel{Hom_{\fu{\cal{C}^{\textrm{op}}}{\cal E}}^{\cal E}}{\longrightarrow} {\cal E}.  
\]
Moreover, the isomorphism above is natural in $F$.
\end{theorem}

\begin{proofs}
One can observe that the internal ${\mathbb{C}}^{\textrm{op}}$-diagram in $\mathbb{E}$ given by the composite
\[
{\mathbb{C}}^{\textrm{op}} \stackrel{ Y^{\textrm{op}} \times \Delta F^{i}  }{\longrightarrow} [{\mathbb C}^{\textrm{op}}, {\cal E}]^{\textrm{op}} \times [{\mathbb C}^{\textrm{op}}, {\cal E}]
\stackrel{Hom_{\fu{\mathbb{C}^{\textrm{op}}}{\cal E}}^{\cal E}}{\longrightarrow} {\cal E}  
\]
is on one hand equal to $F^{i}$ (by an internal version of the Yoneda's lemma) and on the other hand equal to the internalization of the functor $Hom_{\fu{{\cal C}^{\textrm{op}}}{\cal E}}^{\cal E}(\overline{Y}(-),F)$. The verifications are easy and left to the reader.\\
Alternatively, one may proceed as follows.\\
We want to prove that $F(c)\cong Hom_{\fu{{\cal C}^{\textrm{op}}}{\cal E}}^{\cal E}(\overline{Y}(c),F)$, naturally in $c\in \cal C$ (and in $F$). It suffices to observe that we have the following sequence of natural bijections:
\[
\begin{array}{c}
E \longrightarrow Hom_{\fu{{\cal C}^{\textrm{op}}}{\cal E}}^{\cal E}(\overline{Y}(c),F) \\
\hline\\
E^{\ast}\circ \overline{Y}(c) \Longrightarrow E^{\ast}\circ F\\
\hline\\  
\gamma_{{\cal{E}}\slash E}^{\ast} \circ Y(c) \Longrightarrow E^{\ast}\circ F\\
\hline\\
Y(c) \Longrightarrow \Gamma_{{\cal{E}}\slash E} \circ E^{\ast}\circ F\\
\hline\\
\textrm{element of } (\Gamma_{{\cal{E}}\slash E} \circ E^{\ast}\circ F)(c)\\
\hline\\
E \longrightarrow F(c).
\end{array}
\]
\end{proofs}
In the case of flat functors, the theorem specializes to the following result.
\begin{corollary}\label{yrf}
Let $\cal{C}$ be a small category and $\cal{E}$ be a locally small cocomplete topos. Then for every flat functor $\stp{F}{{\cal C}^{\textrm{op}}}{\cal E}$, there is a natural isomorphism of functors
\[
F \cong Hom_{{\bf Flat}({{\cal C}^{\textrm{op}}},{\cal E})}^{\cal E}(\overline{Y}(-),F),
\]
where $\stp{\overline{Y}}{\cal C}{{\bf Flat}({\cal C}^{\textrm{op}},{\cal E})}$ is the functor given by the composite
\[
{\cal{C}} \stackrel{Y}{\longrightarrow} {\bf Flat}({\cal C}^{\textrm{op}}, \Set)
\stackrel{\gamma_{\cal{E}}^{\ast} \circ -}{\longrightarrow} {\bf Flat}({\cal C}^{\textrm{op}}, {\cal E}) 
\]
and $Hom_{{\bf Flat}({\cal C},{\cal E})}^{\cal E}(\overline{Y}(-),F)$ is the functor given by the composite
\[
{\cal{C}}^{\textrm{op}} \stackrel{{\overline{Y}}^{\textrm{op}} \times \Delta F}{\longrightarrow} {{\bf Flat}({\cal C}^{\textrm{op}}, {\cal E})}^{\textrm{op}}\times {\bf Flat}({\cal C}^{\textrm{op}}, {\cal E})
\stackrel{Hom_{{\bf Flat}({\cal{C}}^{\textrm{op}}, {\cal E})}^{\cal E}}{\longrightarrow} {\cal E}.  
\]
\end{corollary}
\begin{proofs}
This immediately follows from the theorem and the remarks in the first section.
\end{proofs}
From now on we will refer to this result as to the Yoneda representation of flat functors.

\section{Representation problems}
In this section we introduce the notion of representation problem in the general context of locally small indexed categories. This concept will lead to a universal characterization of the Yoneda embeddings, which will be employed in the next section to derive a criterion for a theory to be of presheaf type.\\  
\begin{definition}
Let $\cal S$ a cartesian category and $\underline{\mathbb D}$ a locally small $\cal S$-indexed category. A $\cal S$-indexed functor $\stp{F}{{\underline{\mathbb D}}^{\textrm{op}}}{\mathbb{S}}$ is said to be $\cal S$-representable if there exists an object $A \in \cal D$ such that $F$ is isomorphic in $[{\underline{\mathbb D}}^{\textrm{op}}, \mathbb{S}]$  to the composite
\[
{{\underline{\mathbb D}}^{\textrm{op}}} \stackrel{1_{\underline{\mathbb D}} \times \Delta A}{\longrightarrow} {\underline{\mathbb D}}^{\textrm{op}}\times \underline{\mathbb D}
\stackrel{Hom_{{\underline{\mathbb D}}}^{\mathbb S}}{\longrightarrow} {\mathbb{S}}.  
\]
We denote this composite by $Hom_{{\underline{\mathbb D}}}^{\mathbb S}(-,A)$.\\ If $\cal D$ is the underlying category of a locally small $\cal S$-indexed category $\mathbb D$ then we say that a functor $F:{\cal D}^{\textrm{op}}\rightarrow {\cal S}$ is $\cal S$-representable if it is the underlying functor of an indexed functor of the form $Hom_{{\underline{\mathbb D}}}^{\mathbb S}(-,A)$.  
\end{definition}

\begin{definition}
Let $\cal S$ a cartesian category, $\underline{\mathbb D}$ a locally small $\cal S$-indexed category and $\underline{\mathbb K}$ a $\cal S$-indexed full subcategory of $[{\underline{\mathbb D}}^{\textrm{op}}, \mathbb{S}]$. A locally small $\cal S$-indexed category $\underline{\mathbb F}$ together with $\cal S$-indexed functors $\stp{i}{\underline{\mathbb D}}{\underline{\mathbb{F}}}$ and $\stp{r}{\underline{\mathbb K}}{\underline{\mathbb F}}$ is said to be a solution to the $1$-representation problem for $\underline{\mathbb K}$ if $Hom_{{\underline{\mathbb F}}}^{\mathbb S}(-, r(F))\circ i^{\textrm{op}} \cong F$ canonically in $F \in \cal K$. 
\[  
\xymatrix {
{\underline{\mathbb D}}^{\textrm{op}} \ar[d]_{i^{\textrm{op}}} \ar[r]^{F} & \mathbb S \\
{\underline{\mathbb{F}}}^{\textrm{op}} \ar@{-->}[ur]_{Hom_{{\underline{\mathbb{F}}}}^{\mathbb S}(-, r(F))} & }
\]\\  
$(\underline{\mathbb{F}}, \stp{i}{\underline{\mathbb{D}}}{\underline{{\mathbb{F}}}}, \stp{r}{\underline{\mathbb K}}{\underline{\mathbb F})}$ is said to be a solution to the representation problem for $\mathbb K$ if for each $I\in \cal S$ the triple $(\underline{\mathbb{F}}\slash I, \stp{i\slash I}{\underline{\mathbb{D}}\slash I}{\underline{{\mathbb{F}}}\slash I}, \stp{r\slash I}{\underline{\mathbb K}\slash I}{\underline{\mathbb F}\slash I)}$ is a solution to the $1$-representation problem for the $\cal S$$\slash I$-indexed category $\underline{\mathbb K}\slash I$.\\
A solution $(\underline{\mathbb{F}}, \stp{i}{\underline{\mathbb{D}}}{\underline{{\mathbb{F}}}}, \stp{r}{\underline{\mathbb K}}{\underline{\mathbb F})}$ to the representation problem for $\mathbb K$ is said to be universal if for any other solution $(\underline{\mathbb{F'}}, \stp{i'}{\underline{\mathbb{D}}}{\underline{{\mathbb{F'}}}}, \stp{r'}{\underline{\mathbb K}}{\underline{\mathbb F'})}$ to the same problem there exixts a unique (up to canonical isomorphism) $\cal S$-indexed functor $\stp{z}{\underline{\mathbb F}}{\underline{\mathbb F'}}$ such that $z \circ r \cong r'$ and $z \circ i \cong i'$ canonically. Of course, if such a solution exists, it is unique up to canonical isomorphism by the universal property.    
\end{definition} 
     
\begin{proposition}
Let $\cal S$ a cartesian category, $\mathbb{D}$ a locally small $\cal S$-indexed category and $\mathbb K$ a $\cal S$-indexed full subcategory of $[{\underline{\mathbb{D}}}^{\textrm{op}}, \mathbb{S}]$. If $\stp{\mathbb Y}{\mathbb D}{[{\underline{\mathbb{D}}}^{\textrm{op}}, \mathbb{S}]}$ factors as $\stp{\mathbb Y'}{\mathbb D}{\mathbb K}$ through the full embedding $\mathbb K \hookrightarrow [{\underline{\mathbb{D}}}^{\textrm{op}}, \mathbb{S}]$, then the triple $(\underline{\mathbb{K}}, \stp{\mathbb Y'}{\underline{\mathbb{D}}}{\underline{{\mathbb{K}}}}, \stp{1_{\mathbb K}}{\underline{\mathbb K}}{\underline{\mathbb K})}$ is the universal solution to the representation problem for the $\cal S$-indexed category $\mathbb K$. 
\end{proposition}
\begin{proofs}
This is an immediate consequence of the indexed version of the Yoneda lemma.
\end{proofs}
So, if $\mathbb C$ is an internal category in $\cal S$, the embedding $\stp{\mathbb Y}{\mathbb C}{[\mathbb{C}^{\textrm{op}}, \mathbb{S}]}$ can be characterized not only, as it is well known, as the free $\cal S$-cocompletion of $\mathbb C$, but also as the universal solution to the representation problem for the $\cal S$-indexed category $[\mathbb{C}^{\textrm{op}}, \mathbb{S}]$. 

\begin{corollary}\label{corollary}
Let $\mathbb C$ be an internal category in a topos $\cal E$. Then the factorization $\stp{\mathbb Y'}{\mathbb C}{\underline{\bf{Tors}({\mathbb{C}}^{\textrm{op}},\mathbb{E})}}$ of the Yoneda indexed functor $\stp{\mathbb Y}{\mathbb C}{[\mathbb{C}^{\textrm{op}}, \mathbb{E}]}$ through the full embedding $\underline{\bf{Tors}({\mathbb{C}}^{\textrm{op}},\mathbb{E})} \hookrightarrow [{\mathbb{C}}^{\textrm{op}}, \mathbb{E}]$ is the universal solution to the representation problem for the $\cal E$-indexed category $\underline{\bf{Tors}({\mathbb{C}}^{\textrm{op}},\mathbb{E})}$. 
\end{corollary}\qed

\section{Classifying toposes}
As promised, we give a characterization of the (geometric) theories of presheaf type based on the ideas in the last section.\\
We observe that, if $\mathbb T$ is a geometric theory, we can regard it informally as a category $\underline{{\mathbb T}\textrm{-mod}}$ indexed by the (meta)category of Grothendieck toposes via the pseudofunctor $\mathbb T\textrm{-mod}$ (which assigns to every topos $\cal{E}$ the category of $\mathbb T$-models in $\cal{E}$); in particular, for each Grothendieck topos $\cal E$, by ``restricting'' this pseudofunctor to the slices of $\cal E$, we obtain a $\cal E$-indexed category $\underline{{\mathbb T}\textrm{-mod}_{\cal E}}$, which is locally small as a $\cal{E}$-indexed category. Indeed, it is well known that $\mathbb T$ is Morita-equivalent (that is, has the same category of models - up to natural equivalence - into every Grothendieck topos $\cal{E}$ naturally in $\cal{E}$, equivalently has the same classifying topos) to the theory of flat functors on a category $\cal{C}$ which are continuous with respect to a Grothendieck topology $J$ on $\cal{C}$, and the categories of such functors are all full subcategories of the corresponding categories of functors on $\cal{C}$ (cfr. the remarks in the first section).\\    
We recall that a geometric theory $\mathbb T$ is said to be of presheaf type if its classifying topos is a presheaf topos (equivalently, the topos $[{\cal C},\Set]$, where ${\cal C}:=(\textrm{f.p.} {\mathbb T}\textrm{-mod}(\Set))$ is the category of finitely presentable $\mathbb T$-models in $\Set$). If $\mathbb T$ is of presheaf type, then it is Morita-equivalent to the theory of flat functors on the category ${\cal C}^{\textrm{op}}$. Via this equivalence, the Yoneda embedding $\stp{\overline{Y}}{\cal C}{{\bf Flat}({\cal C}^{\textrm{op}},\cal E)}$ corresponds to the embedding of $(\textrm{f.p.} {\mathbb T}\textrm{-mod}(\Set))$ into ${\mathbb T}\textrm{-mod}(\cal E)$ given by the inverse image functor $\gamma_{\cal E}^{\ast}$. Notice that the image in ${\mathbb T}\textrm{-mod}(\cal E)$ of this embedding can be thought as the subcategory of ``constant $\mathbb T$-models which are finitely presentable in $\Set$''.

As we have remarked, for each Grothendieck topos $\cal{E}$ $\underline{{\mathbb T}\textrm{-mod}_{\cal E}}$ is locally small, so it does make sense to ask if $\stp{\gamma^{\ast}_{\cal E}(-)}{\textrm{f.p.} {\mathbb T}\textrm{-mod}(\Set)}{{\mathbb T}\textrm{-mod}(\cal E)}$ (regarded here as a $\cal E$-indexed functor to $\underline{{\mathbb T}\textrm{-mod}_{\cal E}}$) is the universal solution to the representation problem for the $\cal E$-indexed category $\underline{{\bf Flat}((\textrm{f.p.} {\mathbb T}\textrm{-mod}(\Set))^{\textrm{op}}, {\mathbb E})}$. If this holds for every $\cal E$ naturally in $\cal E$ then we may conclude by Corollary \ref{corollary} that $\mathbb T$ is of presheaf type. More concretely, we have the following criterion for a theory to be of presheaf type. 
\begin{theorem}
Let $\mathbb T$ a geometric theory. Then $\mathbb T$ is of presheaf type if and only if for each Grothendieck topos $\cal E$, every flat functor\\ $\stp{F}{(\textrm{f.p.} {\mathbb T}\textrm{-mod}(\Set))^{\textrm{op}}}{\cal E}$\\
can be extended to a $\cal E$-representable along\\ $\stp{{\gamma^{\ast}_{\cal E}(-)}^{\textrm{op}}}{(\textrm{f.p.} {\mathbb T}\textrm{-mod}(\Set))^{\textrm{op}}}{{\mathbb T}\textrm{-mod}(\cal E)^{\textrm{op}}}$\\ and conversely every $\cal E$-representable ${\mathbb T}\textrm{-mod}(\cal E)^{\textrm{op}}\rightarrow \cal E$ arises up to isomorphism in this way, naturally in $F$ and $\cal E$.  
\end{theorem}      
\begin{proofs}
This is immediate from the discussion above.
\end{proofs}

Suppose now you have a geometric theory $\mathbb T$ classified by the topos $[{\cal C}^{\textrm{op}}, \Set]$ and want to understand what the theory $\mathbb T'$ classified by the topos $\Sh ({\cal C}, J )$ (where $J$ is a Grothendieck topology on $\cal C$) looks like, in terms of $\mathbb T$, and without any reference to flat functors. The technique of the Yoneda representation for flat functors provides us with a means for solving this problem. Specifically, we are able to describe in terms of the $\mathbb T$-models and of the Grothendieck topology $J$ the $\mathbb T'$-models in any Grothendieck topos $\cal E$, in other words we are able to identify $\mathbb T'$ up to Morita-equivalence entirely in terms of $\mathbb T$ and of $J$.\\

We denote by $\cal{\check{C}}$ the Cauchy completion of the category $\cal C$. Recall that $\cal{\check{C}}$ can alternatively be characterized as the full subcategory of $\Ind{\cal C}$ consisting of the finitely presentable objects and also as the closure of $\cal C$ under retracts in $\Ind{\cal C}$.\\
It is well known that the functor categories $[{\cal C}^{\textrm{op}}, \Set]$ and $[{\cal{\check{C}}}^{\textrm{op}}, \Set]$ are naturally equivalent. Since $\Sh ({\cal C}, J )$ is a subtopos of $[{\cal{\check{C}}}^{\textrm{op}}, \Set]$, it follows (from the theory of elementary toposes) that there exists a unique Grothendieck topology $\check{J}$ on $\cal{\check{C}}$ such that the toposes $\Sh ({\cal C}, J )$ and $\Sh ({\cal{\check{C}}}, \check{J})$ are naturally equivalent. We describe it explicitly in the theorem below.\\
We adopt the following conventions: if $S$ is a sieve in $\cal C$, we denote by $\overline{S}$ the sieve in $\cal{\check{C}}$ generated by the members of $S$; if $R$ is a sieve in $\cal{\check{C}}$, we denote by $R\cap arr({\cal C})$ the sieve in $\cal C$ formed by the elements of $R$ which are arrows in $\cal C$. Moreover, given an arrow $\stp{g}{d}{c}$ in $\cal C$ and sieves $S$ and $R$ on $c$ respectively in $\cal C$ and $\cal{\check{C}}$, we denote by $g^{\ast}_{\cal C}(S)$ and $g^{\ast}_{\cal{\check{C}}}(R)$ the sieves obtained by pulling back $S$ and $R$ along $g$ respectively in the categories $\cal C$ and $\cal{\check{C}}$.
\begin{theorem}
Let $\cal C$ be a category and $\cal{\check{C}}$ its Cauchy completion.
Given a Grothendieck topology $J$ on $\cal C$, there exists a unique Grothendieck topology $\check{J}$ on $\cal{\check{C}}$ that induces $J$ on $\cal C$, which is defined by:
for each sieve $R$ on $d \in \cal{\check{C}}$, $R \in \check{J}(d)$ if and only if there exists a retract $d \stackrel{i}{\hookrightarrow} a \stackrel{r}{\rightarrow} d$ with $a \in \cal C$ and a sieve $S \in J(a)$ such that $R=i^{\ast}(\overline{S})$.\\
Furhermore, if $d\in {\cal C}$ then $R \in \check{J}(d)$ if and only if there exists a sieve $S$ in ${\cal C}$ on $d$ such that $R=\overline{S}$.\\     
\end{theorem}         
\begin{proofs}
Since the full embedding $\cal C \hookrightarrow \cal{\check{C}}$ is (trivially) dense with respect to every Grothendieck topology on $\cal{\check{C}}$, it follows from the Comparison Lemma (Theorem C2.2.3 in \cite{El2}) and the remarks above that there is at most one Grothendieck topology on $\cal{\check{C}}$ that induces $J$ on $\cal{C}$. Therefore, it will be enough to prove that the coverage $\check{J}$ in the statement of the theorem is a Grothendieck topology that induces $J$ on $\cal C$.\\ This, as well as the second part of the thesis, can be easily proved by using the following easy fact (whose proof is left to the reader):\\
Given an object $c \in {\cal C}$, the assignments $R\rightarrow R\cap arr({\cal C})$ and $S \rightarrow \overline{S}$ are inverse to each other and define a bijection between the set of sieves in $\cal{C}$ on $c$ and the set of sieves in $\cal{\check{C}}$ on $c$. Moreover, these bijections are natural with respect to the operations of pullback of sieves along an arrow in ${\cal C}$.\\ 
By way of example, we provide the details of the proof that $\check{J}$ satisfies the ``stability axiom'' for Grothendieck topologies.\\
Given $R \in \check{J}(d)$ and $\stp{g}{e}{d}$ in $\cal{\check{C}}$, we want to prove that $g^{\ast}(R)\in \check{J}(e)$. 
Since $R \in \check{J}(d)$, there exists a retract $d \stackrel{i}{\hookrightarrow} a \stackrel{r}{\rightarrow} d$ with $a \in \cal C$ and a sieve $S \in J(a)$ such that $R=i^{\ast}(\overline{S})$. There exists a retract $e \stackrel{j}{\hookrightarrow} b \stackrel{z}{\rightarrow} e$ with $b \in \cal C$. Now, $g^{\ast}(R)=g^{\ast}(i^{\ast}(T))=(i\circ g)^{\ast}(T)=((i\circ g\circ z)\circ j)^{\ast}(T)=j^{\ast}((i\circ g \circ z)^{\ast}(T))=j^{\ast}((i\circ g \circ z)_{\cal{\check{C}}}^{\ast}(\overline{S}))=j^{\ast}(\overline{(i\circ g \circ z)^{\ast}_{\cal C}(S)})$. Our thesis then follows at once from the stability axiom for $J$.  
\end{proofs}
\begin{theorem}\label{teoat}
Let $\cal C$ be a category and $J$ a Grothendieck topology on $\cal C$.\\
If $J$ is the trivial topology then $\check{J}$ is the trivial topology.\\
If $J$ is the dense (respectively, the atomic) topology on $\cal C$, then $\check{J}$ is the dense (respectively, the atomic) topology on $\cal{\check{C}}$. 
\end{theorem}
\begin{proofs}
All can be easily proved by using the ``retract technique'' employed in the proof of the previous theorem. We omit the details. 
\end{proofs}

Coming back to our original problem, we have seen that it is natural to replace the topos $\Sh({\cal C},J)$ with $\Sh({\cal{\check{C}}}, \check{J})$. The advantage for us of this replacement is that the category $\cal{\check{C}}$, being Cauchy complete, can be recovered from ${\bf Flat}(\cal{\check{C}}^{\textrm{op}}, \Set)$ as the full subcategory of finitely presentable objects. Hence, if $\mathbb T$ is a theory classified by $[\cal{\check{C}}, \Set]$, then the natural equivalence ${\bf Flat}({\cal{\check{C}}}^{\textrm{op}}, \Set)\simeq {\mathbb T}\textrm{-mod}(\Set)$ restricts to a natural equivalence ${\cal{\check{C}}}\simeq \textrm{f.p.}{\mathbb T}\textrm{-mod}(\Set)$, as in the following diagram:
\[  
\xymatrix {
\cal{\check{C}} \ar[d]_{\overline{Y}} & \textrm{f.p.}{\mathbb T}\textrm{-mod}(\Set) \ar[l]^{\sim}_{\tau} \ar[d]^{i} \\
{\bf Flat}({\cal{\check{C}}}^\textrm{op},\Set) \ar[r]_{\sim} & {\mathbb T}\textrm{-mod}(\Set) } 
\]
Now we want to rewrite the Yoneda representation
\[
F \cong Hom_{\fu{\cal{\check{C}}^{\textrm{op}}}{\cal E}}^{\cal E}(\overline{Y}(-),F),
\]
of a flat functor $F:{\cal{\check{C}}}\rightarrow {\cal E}$ (given by Corollary \ref{yrf}) in terms of $\mathbb T$, regarded here as a $\cal E$-indexed category. We recall that $\underline{{\mathbb T}\textrm{-mod}_{\cal E}}$ is locally small, with $Hom_{{\mathbb T}\textrm{-mod}({\cal E})}^{\cal E}(M,N)$ object of morphisms in $\underline{{\mathbb T}\textrm{-mod}_{\cal E}}$ from $M$ to $N$ belonging to ${\mathbb T}\textrm{-mod}(\cal E)$. The naturality in $\cal E$ of the Morita-equivalence between $\mathbb T$ and the theory of flat functors on ${\cal{\check{C}}}^{\textrm{op}}$ implies the commutativity of the following diagram:
\[  
\xymatrix {
{\bf Flat}({\cal{\check{C}}}^\textrm{op},\Set) \ar[r]^{\sim} \ar[d]_{\gamma^{\ast}_{\cal E}\circ -} & {\mathbb T}\textrm{-mod}(\Set) \ar[d]^{\gamma^{\ast}_{\cal E}(-)} \\
{\bf Flat}({\cal{\check{C}}}^\textrm{op},{\cal E}) \ar[r]_{\sim} & {\mathbb T}\textrm{-mod}(\cal E) } 
\]
From the commutativity of the two diagrams above we deduce the following representation for $F\circ \tau$:
\[
F\circ \tau \cong Hom_{{\mathbb T}\textrm{-mod}(\cal E)}^{\cal E}(\gamma^{\ast}_{\cal E}(i(-)),M_{F}),
\]
where $M_{F}$ is the $\mathbb T$-model in $\cal E$ corresponding to $F\in {\bf Flat}({\cal C}^{\textrm{op}}, {\cal E})$ via the Morita-equivalence.\\
This motivates the following definition.

\begin{definition}\label{deffond}
Let $\cal E$ be a locally small cocomplete topos and $\mathbb T$ a theory of presheaf type. Given a Grothendieck cotopology $J$ on ${\cal C}:=\textrm{f.p.}{\mathbb T}\textrm{-mod}(\Set)$, a model $M\in {\mathbb T}\textrm{-mod}(\cal E)$ is said to be \emph{$J\textrm{-homogeneous}$} if for each cosieve $S\in J(c)$ the family of all the arrows
\[  
Hom_{{\mathbb T}\textrm{-mod}(\cal E)}^{\cal E}(\gamma^{\ast}_{\cal E}(i(f)),M): Hom_{{\mathbb T}\textrm{-mod}(\cal E)}^{\cal E}(\gamma^{\ast}_{\cal E}(i(cod(f))),M) \longrightarrow   Hom_{{\mathbb T}\textrm{-mod}(\cal E)}^{\cal E}(\gamma^{\ast}_{\cal E}(i(c)),M)
\]
for $f\in S$, is epimorphic in $\cal E$.
\end{definition}
\begin{rmk}
\emph{It is clear (from the definition of atomic topology) that if $J$ is the atomic cotopology on ${\cal C}$ then a model $M\in {\mathbb T}\textrm{-mod}(\cal E)$ is $J\textrm{-homogeneous}$ if and only if for each arrow $\stp{f}{c}{d}$ in $\cal C$, the arrow \[  
Hom_{{\mathbb T}\textrm{-mod}(\cal E)}^{\cal E}(\gamma^{\ast}_{\cal E}(i(f)),M): Hom_{{\mathbb T}\textrm{-mod}(\cal E)}^{\cal E}(\gamma^{\ast}_{\cal E}(i(d)),M) \longrightarrow   Hom_{{\mathbb T}\textrm{-mod}(\cal E)}^{\cal E}(\gamma^{\ast}_{\cal E}(i(c)),M)
\] 
is an epimorphism in $\cal E$.\\
In this case we will simply say `homogeneous' instead of `$J$-homogeneous'.} 
\end{rmk}
\flushleft
We observe that $M_{F}$ is $J$-homogeneous if and only if $F\circ \tau$ is $J$-continuous. We thus obtain the following theorem.
\begin{theorem}
Let $({\cal C}, J)$ be a site and $\mathbb T$ a theory classified by the topos $[{\cal C}^{\textrm{op}}, \Set]$. Then the topos $\Sh({\cal C}, J)$ classifies the $\mathbb T$-models which are $\check{J}$-homogeneous; that is, given a geometric theory $\mathbb T'$ together with a full and faithful indexed functor $i:\underline{\mathbb T'\textrm{-mod}}\hookrightarrow \underline{\mathbb T\textrm{-mod}}$, then\\ the $\mathbb T'$-models are identified by $i$ with the $\check{J}$-homogeneous $\mathbb T$-models if and only if\\
\textbullet{ }$\mathbb T'$ is classified by the topos $\Sh({\cal C}, J)$ and\\
\textbullet{ }the embedding $i$ is induced via the universal property of the classifying toposes by the inclusion $\Sh({\cal C}, J)\hookrightarrow [{\cal C}^{\textrm{op}}, \Set]$. \end{theorem}\qed

Specializing the theorem to the case of the atomic topology gives the following result.    
\begin{corollary}\label{homog}
Let $({\cal C}, J)$ be an atomic site and $\mathbb T$ a theory classified by the topos $[{\cal C}^{\textrm{op}}, \Set]$. Then the topos $\Sh({\cal C}, J)$ classifies the homogeneous $\mathbb T$-models.
\end{corollary}
\begin{proofs}
This is immediate from the theorem and Theorem \ref{teoat}.
\end{proofs}

Now we want to rephrase in more explicit terms what it means for a model to be $J$-homogeneous; this will be particularly important for the applications.\\
To this end, we first express the condition that a given family of arrows as in Definition \ref{deffond} is epimorphic as a logical sentence in the internal language of the topos, then we use the Kripke-Joyal semantics to spell out what it means for that sentence to be valid in the topos.\\
Recall that if $\cal E$ is a cocomplete topos and $(f_{i}:C_{i}\rightarrow C\mid i\in I)$ is a family of arrows in it indexed by a set $I$, then this family is epimorphic if and only if the logical formula $(\forall y\in C)(\mathbin{\mathop{\textrm{\huge $\vee$}}\limits_{i\in I}} (\exists x\in C_{i}(f_{i}x=y)))$ holds in $\cal E$. Given a class of generators $\cal G$ for $\cal E$, the validity in $\cal E$ of this sentence is in turn equivalent, by the Kripke-Joyal semantics, to the following statement:\\
for each $E\in \cal G$ and $y:E\rightarrow C$ there exists an epimorphic family $(r_{i}:E_{i}\rightarrow E\mid i\in I)$ and generalized elements $(x_{i}:E_{i}\rightarrow C_{i}\mid i\in I)$ such that $y\circ r_{i}=f_{i}\circ x_{i}$ for each $i\in I$. By applying this to the families of arrows in Definition \ref{deffond} and by recalling that the objects $Hom_{{\mathbb T}\textrm{-mod}(\cal E)}^{\cal E}(\gamma^{\ast}_{\cal E}(i(d)),M)$ are the objects of morphisms from $\gamma^{\ast}_{\cal E}(i(d))$ to $M$ in $\underline{\mathbb T\textrm{-mod}_{\cal E}}$, we obtain the following characterization.  
\begin{theorem}\label{teofond}
Let $\cal E$ be a locally small cocomplete topos with a class of generators $\cal G$ and $\mathbb T$ be a theory of presheaf type. Given a Grothendieck cotopology $J$ on ${\cal C}:=\textrm{f.p.}{\mathbb T}\textrm{-mod}(\Set)$, a model $M\in {\mathbb T}\textrm{-mod}(\cal E)$ is $J\textrm{-homogeneous}$ if and only if for each cosieve $S\in J(c)$, object $E\in \cal G$ and arrow $y:E^{\ast}(\gamma^{\ast}_{\cal E}(i(c)))\rightarrow E^{\ast}(M)$ in ${\mathbb T}\textrm{-mod}({\cal E}\slash E)$ there exists an epimorphic family $(p_{f}:E_{f}\rightarrow E, f\in S)$ and for each arrow $f:c\rightarrow d$ in $S$ an arrow $u_{f}:E_{f}^{\ast}(\gamma^{\ast}_{\cal E}(i(d)))\rightarrow E_{f}^{\ast}(M)$ in ${\mathbb T}\textrm{-mod}({\cal E}\slash E)$ such that $p_{f}^{\ast}(y)=u_{f}\circ E_{f}^{\ast}(\gamma^{\ast}_{\cal E}(i(f)))$.
\end{theorem}\qed
Notice that if $\cal E$ is the topos $\Set$ then by taking as class of generators of $\Set$ the class having as its unique element the singleton $1_{\Set}$ we obtain the following result.
\begin{corollary}\label{cor}
Let $\mathbb T$ be a theory of presheaf type. Given a Grothendieck cotopology $J$ on ${\cal C}:=\textrm{f.p.}{\mathbb T}\textrm{-mod}(\Set)$, a model $M\in {\mathbb T}\textrm{-mod}(\Set)$ is $J\textrm{-homogeneous}$ if and only if for each cosieve $S\in J(c)$ and arrow $y:i(c)\rightarrow M$ in ${\mathbb T}\textrm{-mod}(\Set)$ there exists an arrow $f:c\rightarrow d$ in $S$ and an arrow $u_{f}:i(d)\rightarrow M$ in ${\mathbb T}\textrm{-mod}(\Set)$ such that $y=u_{f}\circ i(f)$.
\end{corollary}\qed 
By specializing the theorem and the corollary to the case of the atomic topology one immediately obtains the following results.
\begin{corollary}\label{atom}
Let $\cal E$ be a locally small cocomplete topos with a class of generators $\cal G$ and $\mathbb T$ be a theory of presheaf type. If ${\cal C}:=\textrm{f.p.}{\mathbb T}\textrm{-mod}(\Set)^{\textrm{op}}$ satisfies the right Ore condition then a model $M\in {\mathbb T}\textrm{-mod}(\cal E)$ is homogeneous if and only if for each arrow $f:c\rightarrow d$ in ${\cal C}^{\textrm{op}}$, object $E\in \cal G$ and arrow $y:E^{\ast}(\gamma^{\ast}_{\cal E}(i(c)))\rightarrow E^{\ast}(M)$ in ${\mathbb T}\textrm{-mod}({\cal E}\slash E)$, there exists an object $E_{f}\in {\cal E}$, an epimorphism $p_{f}:E_{f}\twoheadrightarrow E$ and an arrow $u_{f}:E_{f}^{\ast}(\gamma^{\ast}_{\cal E}(i(d)))\rightarrow E_{f}^{\ast}(M)$ in ${\mathbb T}\textrm{-mod}({\cal E}\slash E)$ such that $p_{f}^{\ast}(y)=u_{f}\circ E_{f}^{\ast}(\gamma^{\ast}_{\cal E}(i(f)))$.
\end{corollary}

\begin{corollary}
Let $\mathbb T$ be a theory of presheaf type. If ${\cal C}:=\textrm{f.p.}{\mathbb T}\textrm{-mod}(\Set)^{\textrm{op}}$ satisfies the right Ore condition then a model $M\in {\mathbb T}\textrm{-mod}(\Set)$ is homogeneous if and only if for each arrow $f:c\rightarrow d$ in $\textrm{f.p.}{\mathbb T}\textrm{-mod}(\Set)$ and arrow $y:i(c)\rightarrow M$ in ${\mathbb T}\textrm{-mod}(\Set)$ there exists an arrow $u_{f}:i(d)\rightarrow M$ in ${\mathbb T}\textrm{-mod}(\Set)$ such that $y=u_{f} \circ i(f)$:
\[  
\xymatrix {
i(c) \ar[d]_{i(f)} \ar[r]^{y} & M \\
i(d) \ar@{-->}[ur]_{u_{f}} & }
\]
\end{corollary}\qed
\begin{rmk}\label{remark}
\emph{We observe that under the hypotheses of Definition \ref{deffond} for each topos $\cal E$ and object $E\in {\cal E}$ there is an isomorphism $E^{\ast}(Hom_{{\mathbb T}\textrm{-mod}({\cal E})}^{\cal E}(\gamma^{\ast}_{\cal E}(i(c)),M))\cong Hom_{{\mathbb T}\textrm{-mod}({\cal E}\slash E)}^{{\cal E}\slash E}(\gamma^{\ast}_{{\cal E}\slash E}(i(c)),E^{\ast}(M))$, which is natural in $c\in \cal C$. Hence, if $M\in {\mathbb T}\textrm{-mod}({\cal E})$ is $J$-homogeneous then $E^{\ast}(M)\in {\mathbb T}\textrm{-mod}({\cal E}\slash E)$ is also $J$-homogeneous. This implies that, while dealing with theories $\mathbb T'$ that one wants to prove to satisfy the conditions of Theorem \ref{teofond}, one can restrict to argue with generalized elements defined on $1$, by the localizing principle. This is illustrated in the following example.}
\end{rmk}
\section{An example}
As an application of Corollaries \ref{homog} and \ref{atom}, we prove that the classifying topos for the theory of dense linearly ordered objects without endpoints is given by atomic topos $\Sh({{\bf Ord}^{\textrm{op}}_{fm}}, J)$, where ${\bf Ord}_{fm}$ is the category of finite ordinals and order-preserving injections between them and $J$ is the atomic cotopology on it.\\  
The theory $\mathbb L'$ of dense linearly ordered objects without endpoints is defined over a one-sorted signature having one relation symbol $\lt$ apart from equality, and has the following axioms:\\  
\[
(((x\lt y)\wedge (y\lt x))\: {\vdash_{x,y}}\: \bot ),
\]
\[
(\top \: \vdash_{x,y}\: ((x=y)\vee (x\lt y)\vee (y \lt x))),
\]
\[
(\top \:\vdash_{\empstg} (\exists x)\top),
\]
\[
((x\lt y)\: \vdash_{x,y} (\exists z)((x\lt z)\wedge (z\lt y))) \textrm{ and}
\]
\[
(x \: \vdash_{x} (\exists y,z)((y\lt x)\wedge (x\lt z))).
\]
The first two axioms give the theory $\mathbb L$ of (decidably) linearly ordered objects; it is well-known that this theory is of presheaf type, hence, being ${\bf Ord}_{fm}$ the category of finitely presentable $\mathbb L'$-models in $\Set$, its classifying topos is equivalent to the functor category $[{\bf Ord}_{fm}, \Set]$. Notice also that the category ${\bf Ord}^{\textrm{op}}_{fm}$ satisfies the right Ore condition, so we can equip it with the atomic topology $J$.\\
A model $M\in \mathbb L\textrm{-mod}(\cal E)$ is given by a pair $(I, R)$ where $I$ is an object of $\cal E$ and $R$ is a relation on $I$ satisfying the diagrammatic forms of the first two axioms above. We will prove that for each topos $\cal E$, a model $M=(I,R)\in \mathbb L\textrm{-mod}(\cal E)$ is homogeneous if and only if it is a model of $\mathbb L'$, that is if $(I,R)$ is non-empty, dense and without endpoints; this will imply (by the corollaries) our thesis.\\
In one direction, let us prove that if $M$ is homogeneous then $(I,R)$ is dense. For each object $E\in {\cal E}$, we denote by $\lt_{E}$ is the order induced by $R$ on $Hom_{\cal E}(E,I)$. By the localizing principle (cfr. Remark \ref{remark}), it is enough to prove that if $x,y:1\rightarrow I$ are two generalized elements of $I$ with $x\lt_{1} y$ then there exists an object $E\in \cal E$, an epimorphism $p:E\twoheadrightarrow 1$ and an arrow $z:E\rightarrow I$ such that $x\comp p\lt_{E} z\lt_{E} y\circ p$. Consider the arrow $f:2\rightarrow 3$ in ${\bf Ord}_{fm}$ defined by $f(0)=0$ and $f(1)=2$; the arrows $x$ and $y$ induce, via the assignment $(0\rightarrow x$, $1\rightarrow y)$ and the universal property of the coproduct $\gamma^{\ast}_{\cal E}(2)$, an arrow $\psi:\gamma^{\ast}_{\cal E}(2)\rightarrow I$ in $\mathbb L\textrm{-mod}(\cal E)$. From the homogeneity of $M$ we obtain the existence of an object $E\in \cal E$, an epimorphism $p:E\twoheadrightarrow 1$ and an arrow $\chi: E^{\ast}(\gamma^{\ast}_{\cal E}(3))\rightarrow E^{\ast}(I)$ in $\mathbb L\textrm{-mod}({\cal E}\slash E)$ such that $\chi\circ E^{\ast}(\gamma^{\ast}_{\cal E}(f))=E^{\ast}(\psi)$. Then the composite arrow 
\[
E\cong E^{\ast}(\gamma^{\ast}_{\cal E}(1))\stackrel{E^{\ast}(\gamma^{\ast}_{\cal E}(u))}{\longrightarrow} E^{\ast}(\gamma^{\ast}_{\cal E}(3))\stackrel{\chi}{\rightarrow} E^{\ast}(I)\stackrel{\pi_{I}}{\rightarrow} I,
\]
where $u:1\rightarrow 3$ is the arrow in ${\bf Ord}_{fm}$ which picks out the element $1\in 3$, gives an arrow $z:E\rightarrow I$ with the required properties. The verifications that $(I,R)$ is non-empty and without endpoints are similar and left to the reader.\\
Conversely, we prove that if $M\in \mathbb L'\textrm{-mod}(\cal E)$ then $M$ is homogeneous. Again, by the localizing principle, this amounts to proving that given an arrow $f:n\rightarrow m$ in ${\bf Ord}_{fm}$ and an arrow $\psi:\gamma^{\ast}_{\cal E}(n)\rightarrow I$ in $\mathbb L\textrm{-mod}(\cal E)$, there exists an object $E\in \cal E$, an epimorphism $p:E\twoheadrightarrow 1$ and an arrow $\chi: E^{\ast}(\gamma^{\ast}_{\cal E}(m))\rightarrow E^{\ast}(I)$ in $\mathbb L\textrm{-mod}({\cal E}\slash E)$ such that $\chi\circ E^{\ast}(\gamma^{\ast}_{\cal E}(f))=E^{\ast}(\psi)$. The arrow $\psi$ can be identified, via the universal property of the coproduct $\gamma^{\ast}_{\cal E}(n)$, with a family $(h_{i}:1\rightarrow I \: \mid i\in n)$ of generalized elements of $I$. To find the required arrow $\chi$, we inductively use the fourth or the fifth axioms to obtain, starting from the $h_{i}$, an object $E\in \cal E$, an epimorphism $p:E\twoheadrightarrow 1$ and $m$ generalized elements $(z_{j}:E\rightarrow I \: \mid  j\in m)$ such that for each $i\in n$ $z_{f(i)}=h_{i}\circ p$ and for each $j,j'\in m$ $((j\lt j')\imp (z_{j}\lt_{E}z_{j'}))$. The family $(z_{j}:E\rightarrow I \: \mid j\in m)$ then gives rise to an arrow $\chi: E^{\ast}(\gamma^{\ast}_{\cal E}(m))\rightarrow E^{\ast}(I)$ in $\mathbb L\textrm{-mod}({\cal E}\slash E)$ with the required property.\qed
\vspace{10 mm}
{\bf Acknowledgements:} I am very grateful to my Ph.D. supervisor Peter Johnstone for his support and encouragement.\\


\newpage
\begin{thebibliography}{10}
\bibitem{El} P. T. Johnstone, \emph{Sketches of an Elephant: a topos theory compendium. Vol.1}, vol. 43 of \emph{Oxford Logic Guides} (Oxford University Press, 2002).
\bibitem{El2} P. T. Johnstone, \emph{Sketches of an Elephant: a topos theory compendium. Vol.2}, vol. 44 of \emph{Oxford Logic Guides} (Oxford University Press, 2002).
\bibitem{MM} S. Mac Lane and I. Moerdijk, \emph{Sheaves in geometry and logic: a first introduction to topos theory} (Springer-Verlag, 1992).
\bibitem{PS} R. Par\'e and D. Schumacher, Abstract families and the Adjoint Functor Theorems, in \emph{Indexed categories and their applications}, Lecture Notes in Math. vol. 661 (Springer-Verlag, 1978), 1-125.
\end{thebibliography}
\end{document}